\def\versiondate{9 Aug. 2009}
\input math.macros
\input Ref.macros

\checkdefinedreferencetrue
\continuousfigurenumberingtrue
\theoremcountingtrue
\sectionnumberstrue
\forwardreferencetrue
\citationgenerationtrue
\nobracketcittrue
\hyperstrue
\initialeqmacro

\input\jobname.key
\bibsty{myapalike}

\def\cp{\tau}   
\def\gh{G}  
\def\gp{\Gamma}  
\def\gpe{\gamma}  
\def\verts{{\ss V}}
\def\vertex{{\ss V}}
\def\edges{{\ss E}}
\def\Tr{{\rm Tr}}  
\def\Det{{\rm Det}} 


\def\asym{{\bf h}}   
\def\ip#1{(\changecomma #1)}
\def\bigip#1{\bigl(\changecomma #1\bigr)}

\def\changecomma#1,{#1,\,}
\def\bigchangecomma#1,{#1,\;}
\def\leftchangecomma#1,{#1,\ }

\def\rtd{\rho}  
\def\bp{o}

\def\alg{{\ss Alg}}  
\def\affalg{\overline{\ss Alg}}  
\def\detalg{{\ss DetAlg}}  

\def\cl#1{\overline #1}  
\def\Dz{{\scr D}_0}  
\def\supp{{\rm supp\/}}  
\def\ER(#1,#2){{\cal R}(#1 \leftrightarrow #2)}   
\def\dom{\succcurlyeq}  
\def\wusf{{\ss WUSF}}

\def\BLPSusf{\ref b.BLPS:usf/, hereinafter referred to as BLPS (2001)%
\def\BLPSusf{BLPS \htmllocref{\bibcode{BLPS:usf}}{(2001)}}}

\def\BLPSgip{\ref b.BLPS:gip/, hereinafter referred to as BLPS (1999)%
\def\BLPSgip{BLPS \htmllocref{\bibcode{BLPS:gip}}{(1999)}}}

\ifproofmode \relax \else\head{To appear in {\it Combin. Probab. Comput.}}
{Version of \versiondate}\fi 
\vglue20pt

\title{Identities and Inequalities for Tree Entropy}

\author{Russell Lyons}

\abstract{The notion of tree entropy was introduced by the author as a
normalized limit of the number of spanning trees in finite graphs, but is
defined on random infinite rooted graphs.
We give some new expressions for tree entropy; one uses Fuglede-Kadison
determinants, while another uses effective resistance. We use the latter to 
prove that tree entropy respects stochastic domination. We also prove that
tree entropy is non-negative in the unweighted case, a special case of
which establishes L\"uck's Determinant Conjecture for Cayley-graph
Laplacians.
We use techniques from the theory of operators affiliated to von Neumann
algebras.
}

\bottomII{Primary 
 05C05, 
60C05. 
Secondary
 05C80. 
}
{Asymptotic complexity, graphs, spanning trees,
determinant, Laplacian, trace.}
{Research partially supported by Microsoft Research and
NSF grants DMS-0406017 and DMS-0705518.}

\bsection{Introduction}{s.intro}

The enumeration of spanning trees in a finite graph is a classical subject
dating to the mid 19th century. Asymptotics began to play a role over 100
years later, in the 1960s.
When a sequence of finite graphs converges in an appropriate, but very
general, sense, \ref
b.Lyons:est/ gave a formula for the limit of the numbers of spanning trees
in that sequence of graphs, when normalized appropriately.
This limit was called the tree entropy of the corresponding limit object,
which was a probability measure on rooted infinite graphs.

This new concept of tree entropy allowed \ref b.Lyons:est/ to give simple
proofs of known limits and inequalities, as well as to resolve an open
question of \ref b.McKay:random/ and to easily calculate new limits.
Here, we give some new expressions for tree entropy, in part correcting
some mistakes in \ref b.Lyons:est/.
Tools we use from the theory of operators affiliated to von Neumann
algebras were not available at the time that \ref b.Lyons:est/ was written.
The new tools also enable us to obtain cleaner results with weaker
hypotheses.
Furthermore, we are able to extend an inequality from \ref b.Lyons:est/
that compares the tree entropies of two probability measures when one
stochastically dominates the other.

The notion of tree entropy extends to weighted graphs, but the case of
unweighted graphs is, of course, particularly interesting.
Using our new representation,
we prove that the tree entropy is non-negative for
unweighted graphs, which is not at all obvious from the definition or from
any of its representations.
In fact, the special case of Cayley graphs establishes L\"uck's Determinant
Conjecture for the graph Laplacian.

We give the details of the results of \ref b.Lyons:est/ referred to above
and then
some background on von Neumann algebras and Fuglede-Kadison
determinants in \ref s.back/. We prove that tree entropy is the logarithm
of a Fuglede-Kadison determinant in \ref t.logdet/. This is used to
represent tree entropy with effective resistances in \ref t.elecent/.
Combined with Rayleigh's monotonicity principle, this representation has
the immediate consequence that stochastic domination implies tree-entropy
domination, \ref t.domin/. This consequence is then combined with
information about wired uniform spanning forests to prove in \ref t.nonneg/
that tree entropy is non-negative for unweighted graphs.

\bsection{Background}{s.back}

In order to define 
the notion of convergence of finite graphs used by \ref b.Lyons:est/ that
we referred to, we first recall the following definitions.
A {\bf rooted graph} $(\gh, \bp)$ is a graph $\gh$ with a distinguished vertex
$\bp$ of $\gh$, called the {\bf root}.
A {\bf rooted isomorphism} of rooted graphs is an isomorphism of the
underlying graphs that takes the root of one to the root of the other.
Given a positive integer $R$, a finite rooted graph $H$, and a probability
distribution $\rtd$ on rooted graphs, let $p(R, H, \rtd)$ denote the
probability that $H$ is rooted isomorphic to the ball of radius $R$ about
the root of a graph chosen with distribution $\rtd$.
For a finite graph $\gh$, let $U(\gh)$ denote the distribution of rooted graphs
obtained by choosing a uniform random vertex of $\gh$ as root of $\gh$.
Suppose that $\Seq{\gh_n}$ is a sequence of finite graphs and that $\rtd$ is a
probability measure on rooted infinite graphs.
We say the {\bf random weak limit} of $\Seq{\gh_n}$ is $\rtd$ if for any
positive integer $R$ and any finite graph $H$, 
we have $\lim_{n \to\infty} p\big(R, H, U(\gh_n)\big) = p(R, H, \rtd)$.
This notion was introduced by \ref b.BS:rdl/.
More generally, if $\gh_n$ are random finite graphs, then
we say the {\bf random weak limit} of $\Seq{\gh_n}$ is $\rtd$ if for any
positive integer $R$, any finite graph $H$, and any $\epsilon > 0$, we have
$\lim_{n \to\infty} \Pleft{\big|p\big(R, H, U(\gh_n)\big) - p(R, H, \rtd)\big| >
\epsilon} = 0$.
Note that only the component of the root matters for convergence to $\rtd$.
Thus, we may and shall assume that $\rtd$ is concentrated on connected graphs.

Recall from \ref b.Lyons:est/ that 
the {\bf tree entropy} of a probability measure $\rtd$ on rooted 
infinite graphs is
$$
\asym(\rtd)
:=
\int \Big(\log \deg_\gh(\bp) - \sum_{k \ge 1} {1 \over k}
p_k(\bp;\gh)\Big)
\,d\rtd(\gh, \bp) 
\,.
\label e.asymdefunwt
$$
One of the main theorems of \ref b.Lyons:est/ was Theorem 3.2, which states
the following. Let $\cp(\gh)$ denote the number of spanning trees of a
graph $\gh$.

\procl t.asym
If $\gh_n$ are finite connected graphs with bounded average degree
whose random weak limit is a probability measure $\rtd$ on infinite rooted
graphs, then 
$$
\lim_{n \to\infty} {1 \over |\verts(\gh_n)|} \log \cp(\gh_n)
=
\asym(\rtd)
\,.
$$
The same limit holds in probability
when $\gh_n$ are random with bounded expected average degree.
\endprocl

In the case of regular graphs $\gh_n$ with girth tending to infinity, the
random weak limit is a rooted regular tree (of the same degree);
with additional hypotheses on $\gh_n$, \ref b.McKay/ proved what amounts to
the same limit as in \ref t.asym/ and asked whether these additional
hypotheses were needed. \ref t.asym/ shows that they are not.

The class of probability measures
$\rtd$ that arise as random weak limits of finite networks is contained in
the class of unimodular $\rtd$, which we now define.
They also include each $\rtd$ that is concentrated on a single Cayley graph
with a fixed root. For more details, see \ref b.AL:urn/.
Since we shall use labeled graphs, i.e., networks, we make a definition
that includes them.

\procl d.unimodular 
Let $\rtd$ be a probability measure on rooted networks.
We call $\rtd$ {\bf unimodular} if 
$$
\int \sum_{x \in \vertex(\gh)} f(\gh, \bp, x) \,d\rtd(\gh, \bp)
=
\int \sum_{x \in \vertex(\gh)} f(\gh, x, \bp) \,d\rtd(\gh, \bp)
$$
for all non-negative Borel functions $f$ on locally
finite connected networks with an ordered pair of distinguished vertices
that is invariant in the sense that for any (non-rooted)
network isomorphism $\gpe$ of $\gh$ and any $x, y \in \vertex(\gh)$, we have
$f(\gpe \gh, \gpe x, \gpe y) = f(\gh, x, y)$,
\endprocl

We need the following finite von Neumann algebra from Section 5 of
\ref b.AL:urn/, to
which we refer for more details.
We also refer to \ref b.Lyons:est/ for more background and motivation.
Suppose that $\rtd$ is a unimodular probability measure on (rooted
isomorphism classes of) rooted (connected) networks.
Consider the Hilbert space $H := \int^\oplus \ell^2\big(\vertex(\gh)\big)
\,d\rtd(\gh, \bp)$, a direct integral.
Let $T : (\gh, \bp) \mapsto T_{\gh, \bp}$ be a measurable assignment of
bounded linear operators $T_{\gh, \bp} : \ell^2\big(\vertex(\gh)\big) \to
\ell^2\big(\vertex(\gh)\big)$
with finite supremum of the norms $\|T_{\gh, \bp}\|$.
Then $T$ induces a bounded linear operator $T := T^\rtd := \int^\oplus T_{\gh,
\bp} \,d\rtd(\gh, \bp)$ on $H$ via
$$
T^\rtd : \int^\oplus f_{\gh, \bp}
\,d\rtd(\gh, \bp) \mapsto \int^\oplus T_{\gh, \bp} f_{\gh, \bp}
\,d\rtd(\gh, \bp)
\,.
$$
The norm $\|T^\rtd\|$ of $T^\rtd$ is the $\rtd$-essential supremum of
$\|T_{\gh, \bp}\|$.
Let $\alg$ be the von Neumann algebra of ($\rtd$-equivalence classes of)
such maps $T$ that are equivariant in
the sense that for all network isomorphisms $\phi : \gh_1 \to \gh_2$, all
$\bp_1, x, y \in \verts(\gh_1)$ and all $\bp_2 \in \verts(\gh_2)$,
we have $(T_{\gh_1, \bp_1} \II x, \II y) = (T_{\gh_2, \bp_2} \phi \II x,
\phi \II y)$.
For $T \in \alg$, we have in particular that $T_{\gh, \bp}$ depends on
$\gh$ but not on the root $\bp$, so we simplify our notation and
write $T_\gh$ in place of $T_{\gh, \bp}$.
Recall that if $T$ is a self-adjoint operator on a Hilbert space
$H$, we write $T \ge 0$ if $\ip{T u, u} \ge 0$ for all $u \in H$.
As shown in Section 5 of \ref b.AL:urn/, the functional
$$
\Tr(T) := \Tr_\rtd(T)
:= \Ebig{(T_\gh \II \bp, \II \bp)} := \int (T_\gh \II \bp, \II \bp)
\,d\rtd(\gh, \bp)
$$
is a trace on $\alg$, which is obviously finite.
Write $\affalg$ for the set of closed densely defined operators affiliated
with $\alg$, i.e., those closed densely defined
operators that commute with all unitary operators
that commute with $\alg$; see, e.g., \ref b.KadRing1/, p.~342.

The only networks we consider will be weighted graphs.
Let $\gh = \big(\vertex(\gh), \edges(\gh), w\big)$ be a graph with a
positive weight function $w : \edges(\gh) \to (0, \infty)$.
For $x \ne y \in \vertex(\gh)$, let $\Delta_\gh(x, y) := -\sum_e w(e)$,
where the sum is over all the edges between $x$ and $y$, and
$\Delta_\gh(x, x) := \sum_e w(e)$,
where the sum is over all non-loop edges incident to $x$.
We assume that $\Delta_\gh(x, x) < \infty$ for all $x$.
An unweighted graph corresponds to $w \equiv 1$, in which case
$\Delta_\gh(x, x)$ is the degree of $x$ in $\gh$ (not counting loops).
The associated network random walk has the transition probability from $x$
to $y$ of $-\Delta_\gh(x, y)/\Delta_\gh(x, x)$; this is simple random walk
in the case of unweighted simple graphs.
Let $p_k(o;\gh)$ be the probability that the network random walk on $\gh$
started at $o$ is again at $o$ after $k$ steps.
The extension from \ref b.Lyons:est/ of \ref e.asymdefunwt/ to weighted
graphs is the following:
the {\bf tree entropy} of a probability measure $\rtd$ on rooted 
weighted infinite graphs is
$$
\asym(\rtd)
:=
\int \Big(\log \Delta_\gh(\bp, \bp) - \sum_{k \ge 1} {1 \over k}
p_k(\bp;\gh)\Big)
\,d\rtd(\gh, \bp) 
\label e.asymdef
$$
whenever this integral converges (possibly to $\pm \infty$).

The {\bf (graph) Laplacian} $\Delta_\gh$ defined in the preceding paragraph
determines an operator 
$$
f \mapsto \Big(x \mapsto \sum_{y \in \vertex} \Delta_\gh(x, y) f(y) \Big)
$$
on functions $f : \vertex(\gh) \to \C$ with finite support.
This operator extends by continuity to a bounded linear operator on
all of $\ell^2\big(\vertex(\gh)\big)$ when $\sup_x \Delta_\gh(x, x) < \infty$.
When 
$$
\rtd\hbox{-}\esssup_{(\gh, \bp)} \sup_{x \in \vertex(\gh)}
\Delta_\gh(x, x) < \infty
\,,
\label e.unifbd
$$
then $(\gh, \bp) \mapsto \Delta_\gh$ defines an operator $\Delta \in
\alg(\rtd)$.
It is self-adjoint and positive semi-definite, i.e., $\Delta \ge 0$.
However, in case we do not have such a uniform bound as \ref e.unifbd/,
we proceed as follows.
Let 
$$
\Dz :=
\{ f \in H \st \all {(\gh, \bp)}  |\supp f_{\gh, \bp}| < \infty \}
\,.
$$
The operator $\Delta$ is defined on the dense subspace $\Dz$, where it is
symmetric.
Let $D$ be the diagonal weighted degree operator on $\Dz$, i.e., $D_\gh(x,
x) := \Delta_\gh(x, x)$ and $D_\gh(x, y) := 0$ for $x \ne y$.
Its closure $\cl D$ is easily seen to be self-adjoint and
affiliated with $\alg$.
Let $P$ be the transition operator for the network random walk, which is
obviously in $\alg$.
Define $\delta := \cl D (I-P)$; since $\cl D \in \affalg$ and $I - P \in \alg$,
it follows that $\delta \in \affalg$.
We claim that $\Delta$ is closeable and
that $\delta = \cl \Delta$.
First, an easy calculation shows that $\delta$ and $\Delta$ agree on $\Dz$,
so that $\delta$ extends $\Delta$.
Since $\delta$ is closed, $\Delta$ is closeable.
Therefore $\cl \Delta \in \affalg$ and, furthermore, is self-adjoint
by Lemma 16.4.1 of \ref b.MvN/ (which is the same as Exercise
6.9.53 of \ref b.KadRing2/).
Since $\cl \Delta \subseteq \delta$, it follows that 
$\cl \Delta = \delta$ by Lemma 16.4.2 of \ref b.MvN/ (or Exercise 6.9.54 of
\ref b.KadRing2/).
From now on, we omit the overlines and write more simply
$D$ and $\Delta$ for their closures, $\cl D$ and $\cl \Delta$.

Let $T \in \affalg$ be a self-adjoint operator 
with spectral resolution $E_T$. We define the Borel measure $\mu_{\rtd, T}$
by
$$
\mu_{\rtd, T}(B) :=
\Tr_\rtd\big(E_T(B)\big)
\label e.specmsr
$$
for Borel subsets $B \subseteq \R$.
We extend the trace by defining 
$$
\Tr_\rtd(T) := 
\int_0^\infty \lambda \,d\mu_{\rtd, T}(\lambda) 
$$
for positive operators $T \in \affalg$ and then by linearity to all of
$\affalg$ when it makes sense.
Write $|T| := \sqrt{T^* T}$.

As in \ref b.HaagSch/ (though with different notation),
write $\detalg$ for the set of $T \in \affalg$ for
which 
$$
\Tr_\rtd(\log^+ |T|)
=
\int_0^\infty \log^+ \lambda \,d\mu_{\rtd, |T|}(\lambda) < \infty
\,.
$$
(The equality is justified by the functional calculus; see Theorem 5.6.26
of \ref b.KadRing1/.)
For $T \in \detalg$, we define its {\bf Fuglede-Kadison determinant} by 
$$
\Det (T)
:=
\Det_\rtd(T) 
:=
\exp \int_0^\infty \log \lambda \,d\mu_{\rtd, |T|}(\lambda)
\in \CO{0, \infty}
\,.
\label e.FKdet
$$
For example, for the diagonal weighted degree operator, $D$, its
Fuglede-Kadison determinant is the geometric-mean weighted degree of the
root: 
$$
\Det_\rtd D 
=
\exp \int \log D_\gh(\bp, \bp) \,d\rtd(\gh, \bp)
\label e.detD
$$
provided this is $< \infty$;
this can be seen either
from the definition by using the fact that $\mu_{\rtd, D}$
is the law of $D_\gh(\bp, \bp)$, or
alternatively
by truncation of $D$ and Fubini's theorem.

\bsection{Tree Entropy}{s.ent}

We now give two new representations of tree entropy and two consequences. The
first representation is as the logarithm of a Fuglede-Kadison determinant.

\procl t.logdet
If $\rtd$ is a unimodular probability measure on rooted weighted
connected infinite graphs with 
$$
\int \log D_\gh(\bp, \bp) \,d\rtd(\gh, \bp) \in \CO{-\infty, \infty}
\,,
\label e.logfinite
$$
then 
$$
\asym(\rtd)
=
\log \Det_{\rtd} \Delta
\in \CO{-\infty, \infty}
\,.
\label e.logdet
$$
\endprocl

\proof
The hypothesis is equivalent to $D \in \detalg$. 
Since $I - P \in \alg \subseteq \detalg$, it follows that $\Delta = D (I-P)
\in \detalg$ with 
$$
\Det\, \Delta = \Det\, D \cdot \Det (I-P)
\label e.detparts
$$
by Proposition 2.5 of \ref b.HaagSch/ (which extends the fundamental
theorem of \ref b.FugKad/ to unbounded operators, as well as to
non-invertible operators).

Since $\|P\| \le 1$,
we have for $0 < c < 1$ that $\log |I - c P| \le (\log 2) I$.
Also, $|I - c P| \to |I - P|$ in the strong operator topology as $c
\uparrow 1$, whence $\log |I - c P| \to \log |I - P|$ in the measure
topology (for its definition, see \ref b.FackKos/, \S 1.5).
Thus,
$$
\Det (I-P) = \lim_{c \uparrow 1} \Det (I - c P)
$$
by the Monotone Convergence Theorem; see, e.g., \ref b.FackKos/, Theorem
3.5(ii).
On the other hand, for $0 < c < 1$,
$$
\log \Det (I - c P)
=
\Re \Tr \log (I - c P)
$$
by Theorem 1 ($2^o$) of
\ref b.FugKad/ (or Theorem I.6.10 of \ref b.Dixmier/) and 
$$
\log (I - c P)
=
-\sum_{k \ge 1} c^k P^k/k
$$
(in the norm topology).
Therefore,
$$
\log \Det (I - c P)
=
-\sum_{k \ge 1} \Re \Tr_\rtd c^k P^k/k
=
-\sum_{k \ge 1} \Tr_\rtd c^k P^k/k
\,,
$$
whose limit as $c \uparrow 1$ is
$$
-\sum_{k \ge 1} \Tr_\rtd P^k/k
=
\int - \sum_{k \ge 1} {1 \over k}
p_k(\bp;\gh)
\,d\rtd(\gh, \bp) 
\label e.trpart
$$
by the Monotone Convergence Theorem.
Comparing \ref e.asymdef/ with equations \ref e.detparts/, \ref e.detD/,
and \ref e.trpart/, we deduce the equality in \ref e.logdet/.
\Qed

\procl r.mistake
The version of this theorem given in \ref b.Lyons:est/ was incorrect even
in the case of unweighted graphs, except when the degrees were bounded.
For example, in the notation used there, whenever the degrees are
unbounded, one gets $\Delta_{G_M}(\bp, \bp) = 0$ with positive probability,
which means that $\Det_\rtd(\Delta_{G_M}) = 0$.
However, unbounded-degree graphs are quite natural, arising, for example,
as limits of random finite graphs.
In addition to that mistake, stronger hypotheses were assumed,
which we now see to be superfluous, and the conclusion was less appealing,
being expressed as a double limit.
\endprocl

An example of a unimodular probability measure $\rtd$ satisfying not only
\ref e.logfinite/, but even the stronger
$$
\int |\log D_\gh(\bp, \bp)| \,d\rtd(\gh, \bp) < \infty
\,,
\label e.abslogfinite
$$
yet with $\asym(\rtd) = -\infty$ is the following.
We work on the nearest-neighbor graph of the integers, $\Z$,
rooted at 0.
Define the weight to be 1 of every edge of the form $(2n, 2n+1)$ for $n \in
\Z$.
Let $X$ be an integer-valued random variable such that 
$\P[X \ge m] = 1/\sqrt m$ for $m \ge 1$.
Let $X_n$ be i.i.d.\ copies of $X$ for $n \in \Z$ and let the weight be
$e^{-X_n}$ of the edge $(2n-1, 2n)$.
Define $\rtd$ to be the resulting measure on rooted weighted graphs.
(In fact, $\rtd$ is defined on rooted isomorphism classes of networks, so
that one does not notice the difference between ``even" and ``odd" edges.)
By Theorem 3.2 of \ref b.AL:urn/, $\rtd$ is unimodular.
Since
$$
\int |\log D_\Z(0, 0)| \,d\rtd 
=
\E[\log (1+e^{-X})]
\,,
$$
\ref e.abslogfinite/ is clearly satisfied.
On the other hand, it is easy to see that there are constants $c_1, c_2 >
0$ such that $p_{2k}(0; \Z) \ge c_1$ for $1 \le k \le \exp \min\{X_0, X_1\}$,
whence
$$
\sum_{k \ge 1} {1 \over k} p_k(0; \Z)
\ge
c_2 \min\{X_0, X_1\}
\,.
$$
Therefore 
$$
\int \sum_{k \ge 1} {1 \over k} p_k(0; \Z) \,d\rtd
\ge
c_2 \E[\min\{X_0, X_1\}]
=
c_2 \sum_{m \ge 1} \P[X_0 \ge m]^2
=
\infty
\,.
$$

A small, but significant,
extension of Theorem 4.2 of \ref b.Lyons:est/ is the following.
Let $(G_1, o_1, w_1)$ and $(G_2, o_2, w_2)$ be two rooted weighted graphs.
Say that $(G_1, o_1, w_1)$ {\bf dominates} $(G_2, o_2, w_2)$ if there is a
graph isomorphism $\phi$ from $G_2$ to a subgraph of $G_1$ that takes
$o_2$ to $o_1$ and such that for
all $e \in \edges(G_2)$, we have $w_2(e) \le w_1\big(\phi(e)\big)$.
This notion is a partial order on rooted weighted graphs and we use the usual
notion of stochastic domination that corresponds to it.
That is,
if $\rtd_1$ and $\rtd_2$ are two probability measures on rooted weighted
graphs, 
say that $\rtd_1$ {\bf stochastically dominates} $\rtd_2$ if there exists a
probability measure $\nu$ on pairs $\big((G_1, o_1, w_1), (G_2, o_2,
w_2)\big)$ 
such that the $\nu$-law of $(G_i, o_i, w_i)$ is $\rtd_i$ for $i = 1, 2$
and $(G_1, o_1, w_1)$ dominates $(G_2, o_2, w_2)$
$\nu$-a.s.

\procl t.domin
If $\rtd_1 \ne \rtd_2$ are unimodular probability measures on rooted weighted
connected infinite graphs that both satisfy \ref e.logfinite/
and $\rtd_1$ stochastically dominates
$\rtd_2$, then $\asym(\rtd_1) > \asym(\rtd_2)$.
\endprocl

The proof of the corresponding result, Theorem 4.2, in \ref b.Lyons:est/
was in fact not complete. We give a more direct proof here based on a
different approach.
In addition, Theorem 4.2 of \ref b.Lyons:est/
assumed \ref e.abslogfinite/ in place of
our hypothesis \ref e.logfinite/ and also assumed a further bound.

The significance of our extension is that Theorem 4.2 of \ref b.Lyons:est/
required the two probability measures $\rtd_i$ to be coupled on the {\it
same\/} graphs, differing only in their edge weights. This makes it
impossible to handle naturally occurring stochastic domination situations,
such as those occurring for limits of random finite graphs.
Thus, the present result can answer a question of \ref b.Lyons:est/
concerning the giant component in the Erd\H{o}s-R\'enyi model of random
graphs, provided one can show stochastic domination of Poisson-Galton-Watson
measures conditioned on survival. Indeed, this domination was proved by
\ref b.LPS:GD/.

To prove \ref t.domin/, we rely on an entirely new representation of tree
entropy.
Given a network $\gh$, one of its vertices $x$, and a positive number $s$,
let $R(\gh, x, s)$ be the effective resistance between $x$ and infinity
in the network $\gh^s$ formed from $\gh$ by adding
an edge of conductance $s$ between every vertex and infinity, where
$\infty$ is also a vertex of $\gh^s$.
To be more precise, consider an exhaustion of $\gh$ by finite subnetworks
$\gh_n$.
Let $H_n$ be the network formed from $\gh$ by identifying all vertices
outside $\gh_n$ to a single vertex $z_n$ and then adding an edge of
conductance $s$ between each vertex of $\gh_n$ and $z_n$.
For large enough $n$, we have that $x \in \vertex(\gh_n)$, so that we may
define the effective resistance $\ER(x, z_n; H_n)$ between $x$ and $z_n$ in
$H_n$. These effective resistances have a limit, which we are calling
$R(\gh, x, s)$.

Our second representation of tree entropy is in terms of electrical
resistance.

\procl t.elecent
If $\rtd$ is a unimodular probability measure on rooted weighted
infinite graphs that satisfies \ref e.logfinite/,
then 
$$
\asym(\rtd)
=
\int_0^\infty \left({s \over 1 + s^2} - \int R(\gh, \bp, s)
\,d\rtd(\gh, \bp)
\right)
\,d s
\,.
\label e.trlog
$$
\endprocl

\procl r.altrep
Although one might ask from comparing \ref e.asymdef/ and \ref e.trlog/
whether for every network $(\gh, \bp)$, one has
$$
\log D_\gh(\bp, \bp) - \sum_{k \ge 1} {1 \over k} p_k(\bp;\gh)
\buildrel ? \over =
\int_0^\infty \left({s \over 1 + s^2} - R(\gh, \bp, s) \right) ds
\,,
\label e.altrep
$$
this is not true. 
Thus, \ref t.elecent/ depends crucially on the assumption that $\rtd$ is
unimodular.
One can show, however, that \ref e.altrep/ does hold for every regular
graph $\gh$: If $d := D_\gh(\bp, \bp)$, then one can show that
$R(\gh, \bp, s)$ equals the expected number of visits to $\bp$ divided by
$d+s$, which equals $\sum_{k \ge 0} p_k(\bp; \gh) d^k/(d+s)^{k+1}$. This
gives that $\int_0^\infty \big(R(\gh, \bp, s) - 1/(d+s)\big) d s = \sum_{k
\ge 1} p_k(\bp; \gh)/k$. Combining this with \ref e.logident/ below gives
the result.
\endprocl

\procl r.fixed
One might also ask whether tree entropy increases under stochastic
domination regardless of the unimodularity of $\rtd$.
This is not the case, however. 
For example, consider $\rtd_1$ to be the measure concentrated on the fixed
graph where the root has degree 1, its neighbor has degree 2, and the
neighbor of the root's neighbor has attached a tree of very large degree.
Let $\rtd_2$ be the measure concentrated on the same graph to which has
been adjoined a loop at the root.
Then a straightforward calculation shows that $\asym(\rtd_1) >
\asym(\rtd_2)$, even though $\rtd_2 \dom \rtd_1$.
\endprocl

\proofof t.elecent
For $\lambda > 0$, a well-known identity states that
$$
\log \lambda
=
\int_0^\infty \left( { s \over 1 + s^2} - {1 \over \lambda + s} \right) ds
\,.
\label e.logident
$$
Also, we have the lesser-known identity
$$
{1 \over 2} \log (1 + \lambda^2)
=
\int_0^\infty \left( { s \over 1 + s^2} - {1 \over \lambda + s} \right)^+ ds
\,.
\label e.possidedom
$$

Since $\Delta \ge 0$, the fact that $\asym(\rtd) < \infty$ (by \ref
t.logdet/) implies that
$$
\int_0^\infty \log (1 + \lambda^2) \,d\mu_{\rtd, \Delta}(\lambda) < \infty
\label e.loglambdafinite
$$
by \ref e.logdet/ and \ref e.FKdet/.

For $s > 0$, note that $(\Delta + s I)^{-1} \in \alg$ since $\Delta \ge 0$
and define $v_s := \big(\Delta_{\gh} + s I\big)^{-1} \II \bp$ on
$\vertex(\gh)$. 
We claim that
$\bigip{v_s, \II \bp} = v_s(\bp) = R(\gh, \bp, s)$.
Indeed, the invertibility of $\Delta + s I$ tells us that $v_s$ is the
unique function on $\vertex(\gh)$ that satisfies $(\Delta + s I) v_s = \II
\bp$.
Since one such function is the limit of the voltage functions $v_{s, n}$
corresponding to the unit current flows on $H_n$
from $\bp$ to $z_n$, it follows that $v_s = \lim_{n \to\infty} v_{s, n}$.
Since $v_{s, n}(\bp) = \ER(\bp, z_n; H_n)$, we obtain the claim.
Hence 
$$
\Tr_\rtd \big((\Delta + s I)^{-1}\big)
=
\int R(\gh, \bp, s) \,d\rtd(\gh, \bp)
\,.
\label e.inv
$$
On the other hand,
$$
(\Delta + s I)^{-1}
=
\int_0^\infty (\lambda + s)^{-1} dE_{\Delta}(\lambda)
\,,
$$
so that
$$
\Tr_\rtd \big((\Delta + s I)^{-1}\big)
=
\int_0^\infty (\lambda + s)^{-1} d\mu_{\rtd, \Delta}(\lambda)
\,.
\label e.dps
$$

Therefore, we have 
$$\eqaln{
\asym(\rtd)
&=
\int_0^\infty \log \lambda \,d\mu_{\rtd, \Delta}(\lambda)
=
\int_0^\infty 
\int_0^\infty \left( { s \over 1 + s^2} - {1 \over \lambda + s} \right) ds
\,d\mu_{\rtd, \Delta}(\lambda)
\cr&=
\int_0^\infty 
\int_0^\infty \left( { s \over 1 + s^2} - {1 \over \lambda + s} \right)
d\mu_{\rtd, \Delta}(\lambda)
\,ds
\cr&=
\int_0^\infty 
\left( { s \over 1 + s^2} - \Tr_\rtd (\Delta + s I)^{-1} \right)
d s
\cr&=
\int_0^\infty 
\left( { s \over 1 + s^2} - \int R(\gh, \bp, s) \,d\rtd(\gh, \bp)
\right)
d s
\,;
}$$
we have used \ref e.logdet/ and \ref e.FKdet/
in the first equality; \ref e.logident/ in the
second; \ref e.possidedom/, \ref e.loglambdafinite/,
and Fubini-Tonelli's Theorem in
the third; \ref e.specmsr/ and \ref e.dps/
in the fourth; and \ref e.inv/ in the fifth.
\Qed

\ref t.domin/ follows immediately by Rayleigh's monotonicity principle.
Indeed, that principle gives us that when $(\gh_1, w_1, \bp_1)$ dominates
$(\gh_2, w_2, \bp_2)$, then 
$$
R(\gh_1, \bp_1, s) \le R(\gh_2, \bp_2, s)
$$
for all $s > 0$, where the edge conductances are understood but not notated
in this inequality.

\procl t.nonneg
If $\rtd$ is a unimodular probability measure on rooted infinite (unweighted)
graphs that satisfies \ref e.logfinite/,
then $\asym(\rtd) \ge 0$, with equality
iff $\int D_\gh(\bp, \bp) \,d\rtd(\gh, \bp) = 2$
iff $\rtd$-a.s.\ $\gh$ is a locally finite tree with 1 or 2 ends.
\endprocl

\proof
By Proposition 7.1 of \ref b.AL:urn/, the root in
the wired uniform spanning forest 
of $\rtd$, denoted $\wusf(\rtd)$, has expected degree 2, whence,
by Theorem 6.2 of \ref b.AL:urn/, the unimodular probability measure 
$\wusf(\rtd)$ is concentrated on trees with at most 2 ends.
This implies that $\wusf(\rtd)$ is amenable by Corollary 8.9 of \ref
b.AL:urn/, whence is the random weak limit of finite trees.
Of course, finite trees have average degree less than 2.
By Theorem 3.2 of \ref b.Lyons:est/,
this means that $\asym\big(\wusf(\rtd)\big) = 0$.
Since $\rtd$ clearly stochastically dominates $\wusf(\rtd)$, it follows by
\ref t.domin/ that $\asym(\rtd) \ge 0$.
The equality condition also follows from \ref t.domin/ and the above
argument, combined with Theorem 6.2 of \ref b.AL:urn/ again.
\Qed

\procl r.mistake2
Proposition 4.3 and Theorem 4.4
in \ref b.Lyons:est/ stated the same results as \ref t.nonneg/, though with
an hypothesis far stronger than \ref e.logfinite/.
However, the proofs relied on a result in a
preliminary version of \ref b.AL:urn/ whose proof was incorrect.
\endprocl

In the special case that $\rtd$ is concentrated on a fixed Cayley graph
$\gh$, then \ref t.nonneg/ says that $\Det\, \Delta_\gh \ge 1$.
This establishes a special case of L\"uck's Determinant Conjecture, which
says that for every group $\gp$ and
for every positive self-adjoint finite matrix over the group ring $\Z\gp$, its
Fuglede-Kadison determinant is at least 1; see, e.g., \ref b.Elek:hyper/.

A consequence of \ref t.domin/ is that the set of measures of fixed tree
entropy and satisfying \ref e.logfinite/
form an anti-chain (no two are comparable in the stochastic domination
order). 
In the special case of tree entropy 0, if we combine this with
\ref t.nonneg/, then we obtain 
that the measures on trees with at most 2 ends and satisfying \ref
e.logfinite/ form an anti-chain: 

\procl c.antichain
If $\rtd_1$ and $\rtd_2$ are unimodular probability measures on rooted 
unweighted infinite trees with at most two ends,
both measures satisfy \ref e.logfinite/,
and $\rtd_1$ stochastically dominates
$\rtd_2$, then $\rtd_1 = \rtd_2$.
\endprocl

\medbreak
\noindent {\bf Acknowledgements.}\enspace I am grateful to Hari Bercovici
and Oded Schramm for conversations.

\def\noop#1{\relax}
\input \jobname.bbl

\filbreak
\begingroup
\eightpoint\sc
\parindent=0pt\baselineskip=10pt

Department of Mathematics,
Indiana University,
Bloomington, IN 47405-5701
\emailwww{rdlyons@indiana.edu}
{http://mypage.iu.edu/\string~rdlyons/}

\endgroup

\bye